\documentclass[11pt,longtable,letterpaper,oneside,reqno]{amsart}
\usepackage{amsmath,amsfonts,amsthm,amssymb,amscd,mathtools,stmaryrd,url}
\usepackage{bbm}
\usepackage{epsfig}
\usepackage{graphicx}
\usepackage{latexsym}
\usepackage{longtable}
\usepackage{float}
\usepackage{hyperref}
\usepackage{color}

\DeclareFontFamily{OT1}{rsfs}{}
\DeclareFontShape{OT1}{rsfs}{n}{it}{<-> rsfs10}{}
\DeclareMathAlphabet{\mathscr}{OT1}{rsfs}{n}{it}

\newtheorem{prop}{Proposition}[section]
\newtheorem{thm}[prop]{Theorem}

\newtheorem{lem}[prop]{Lemma}

\newtheorem {defn }{Definition}
\newtheorem{conjec}[prop]{Conjecture}

\numberwithin{equation}{section}
\usepackage[margin=1in]{geometry}
\usepackage{lipsum}

\begin{document}
\title{On near-perfect numbers of special forms}
\author[E. Hasanalizade]{Elchin Hasanalizade}
\address{Department of Mathematics and Computer Science, University of Lethbridge, 4401 University Drive, Lethbridge, Alberta, T1K 3M4 Canada}
\email{e.hasanalizade@uleth.ca}
\begin{abstract}
In this paper we discuss near-perfect numbers of various forms. In particular, we study the existence of near-perfect numbers in the Fibonacci and Lucas sequences, near-perfect values taken by integer polynomials and repdigit near-perfect numbers.
\end{abstract}

\date{\today}
\maketitle

\section{Introduction}

Let $\sigma(n)$ and $\omega(n)$ denote the sum of the positive divisors of $n$ and the number of distinct prime factors of $n$, respectively. A natural number $n$ is {\it perfect} if $\sigma(n)=2n$. More generally, given a fixed integer $k$, the number $n$ is called {\it multiperfect} or {\it $k$-fold perfect} if $\sigma(n)=kn$. The famous Euclid-Euler theorem asserts that an even number is perfect if and only if it has the form $2^{p-1}(2^p-1)$, where both $p$ and $2^p-1$ are primes. It is not known if there are odd perfect numbers.

In 2012, Pollack and Shevelev \cite{PS} introduced the concept of near-perfect numbers. A positive integer $n$ is {\it near-perfect} with redundant divisor $d$ if $d$ is a proper divisor of $n$ and $\sigma(n)=2n+d$. They constructed the following three types of even near-perfect numbers. 

{\bf Type A.} $n=2^{p-1}(2^p-1)^2$ where both $p$ and $2^p-1$ are primes and $2^p-1$ is the redundant divisor.

{\bf Type B.} $n=2^{2p-1}(2^p-1)$ where both $p$ and $2^p-1$ are primes and $2^p(2^p-1)$ is the redundant divisor.

{\bf Type C.} $n=2^{t-1}(2^t-2^k-1)$, $t\ge k+1$ where $2^t-2^k-1$ is prime and  $2^k$ is the redundant divisor.

In 2013, Chen and Ren \cite{RC} proved that all near perfect numbers $n$ with $\omega(n)=2$ are of types $A$, $B$ and $C$ together with 40. It is an open problem to classify all even near-perfect numbers. On the other hand, from the definition it is easy to see that all odd near-perfect numbers are squares. Li et al. \cite{TRL} showed that there is no odd near-perfect number $n$ with $\omega(n)=3$ and Feng et al. \cite{TMF} proved that the only odd near-perfect number $n$ with $\omega(n)=4$ is $173369889=3^4\cdot7^2\cdot11^2\cdot19^2$. Thus for any other odd near-perfect number $n$ if it exists we have $\omega(n)\ge5$.

There are several papers discussing perfect and multiperfect numbers of various forms. For example, Luca \cite{Luca} proved that there are no perfect Fibonacci or Lucas numbers, while Broughan et al. \cite{BGLLHT} showed that no Fibonacci number (larger than 1) is multiperfect. Assuming the $ABC$-conjecture, Klurman \cite{Klurman} proved that any integer polynomial of degree $\ge3$ without repeated factors can take only finitely many perfect values. Pollack \cite{Pollack} studied perfect numbers with identical digits in base $g$, $g\ge2$. He found that in each base $g$ there are only finitely many examples and that when $g=10$, the only example is $6$. Later, Luca and Pollack \cite{LP} established the same results for multiperfect numbers with identical digits.

In this short note we study near-perfect numbers in the Fibonacci and Lucas sequences, near-perfect values taken by the integer polynomials and near-perfect numbers with identical digits. Recall that the {\it Fibonacci sequence} $(F_n)_{n\ge0}$ is given by $F_0=0$, $F_1=1$ and $F_{n+2}=F_{n+1}+F_n$  for all $n\ge0$ and the {\it Lucas sequence} $(L_n)_{n\ge0}$ is given by $L_0=2$, $L_1=1$ and $L_{n+2}=L_{n+1}+L_n$  for all $n\ge0$. A natural number is called a {\it repdigit in base g} if all of the digits in its base $g$ expansion are equal.

Here we prove the following results

\begin{thm}

a) There are no odd near-perfect Fibonacci or Lucas numbers;

b) There are no near-perfect Fibonacci numbers $F_n$ with $\omega(F_n)\le3$;

c) The only near-perfect Lucas number $L_n$ with two distinct prime factors is $L_6=18$.
\end{thm}

\begin{thm}
\label{Thm2}
Suppose $P(x)\in\mathbb{Z}[x]$ with $\text{deg}P(x)\ge3$ has no repeated factors. Then there are only finitely many $n$ such that $P(n)$ is an odd near-perfect number. Furthermore, if we assume that the $ABC$-conjecture holds, then $P(n)$ takes only finitely many near-perfect values with two distinct prime factors.

\end{thm}

\begin{thm}
\label{Thm3}
Let $2\le g\le10$. Then:
 
a) There are only finitely many repdigits in base $g$ which are near-pefect and have two distinct prime factors. All such numbers are strictly less than $\frac{g^3-1}{g-1}$. In particular, when $g=10$, the only repdigit near-perfect number with two distinct prime divisors is $88$;

b) There are no odd near-perfect repdigits in base $g$.
\end{thm}

\section{Preliminary results}

In this section we collect several auxilary results. We begin with the famous and remarkable theorem of Bugeaud et al. \cite{BMS} about perfect powers in the Fibonacci and Lucas sequences.

\begin{thm}(Bugeaud-Mignotte-Siksek)
\label{BMS}
The only perfect powers among the Fibonacci numbers are $F_0=0$, $F_1=F_2=1$, $F_6=8$, and $F_{12}=144$. For the Lucas numbers, the only perfect powers are $L_1=1$ and $L_3=4$.
\end{thm}

In \cite{Pongsriiam}, Pongsriiam gave the description of the Fibonacci numbers satisfying $\omega(F_n)=3$. We state his results in the following theorems

\begin{thm}
\label{P1}
The only solutions to the equation $\omega(F_n)=3$ are given by 
\begin{gather*}
n=16, 18, \ \text{or} \ 2p \ \text{for some prime} \ p\ge19,\\ 
n=p, p^2, p^3 \ \text{for some prime} \ p\ge5,\\
n=pq \ \text{for some distinct primes} \ p, q\ge3  
\end{gather*}
\end{thm}

\begin{thm}
\label{P2}
Assume that $\omega(F_n)=3$ and $n=p_1p_2$ where $p_1<p_2$ are odd primes. Then $F_{p_1}=q_1$, $F_{p_2}=q_2$ and $F_n=q^{a_1}_1q_2q^{a_3}_3$ where $q_1, q_2, q_3$ are distinct primes, $q_3$ is a primitive divisor of $F_n$ (i.e. a prime divisor which does not divide any $F_m$ for $0<m<n$), $a_3\ge1$ and $a_1\in\{1,2\}$. Furthemore $a_1=2$ if and only if $q_1=p_2$. 
\end{thm}

Let us also recall the $ABC$-conjecture. For $n\in\mathbb{Z}\setminus\{0\}$ the {\it radical} of $n$ is defined by $\text{rad}(n)=\prod_{p|n}p$.

\begin{conjec}($ABC$-conjecture)
For each $\epsilon>0$, there exists $M_{\epsilon}>0$ such that whenever $a,b,c\in\mathbb{Z}\setminus\{0\}$ satisfy the conditions
\begin{align*}
\text{gcd}(a,b,c)=1 \ \text{and} \ a+b=c
\end{align*}
then 
\begin{align*}
\text{max}\{|a|, |b|, |c|\}\leq M_{\epsilon}\text{rad}(abc)^{1+\epsilon}
\end{align*}
\end{conjec}

The next lemma is important for the proof of Theorem \ref{Thm2}.

\begin{lem}(\cite[Corollary 2.4]{Klurman})
Assume that the $ABC$-conjecture is true. Suppose that $f(x)\in\mathbb{Z}[x]$ is non-constant and has no repeated roots. Fix $\epsilon>0$, then 
\begin{align}
\label{Eq2.1}
\prod_{p|f(m)}p\gg |m|^{\text{deg} f-1-\epsilon}
\end{align}
\end{lem}

We also need the finiteness result for the solutions of the hyperelliptic equation
\begin{thm}(Baker \cite{Baker})
\label{B}
All solutions in integers $x,y$ of the diophantine equation
\begin{align*}
y^2=a_0x^n+a_1x^{n-1}+\ldots+a_n
\end{align*}
where $n\ge3$, $a_0\ne0$, $a_1,\ldots,a_n$ denote rational integers and where the polynomial on the right of possesses at least three simple zeros satisfy 
\begin{align*}
\text{max}(|x|,|y|)<\text{exp}\text{exp}\text{exp}\{(n^{10n}\mathcal{H})^{n^2}\}
\end{align*}
where $\mathcal{H}=\text{max}_{0\le j\le n}|a_j|$.
\end{thm}

The next two theorems characterize those perfect powers which are also repdigits.

\begin{thm}(Bugeaud-Mignotte \cite{BM})
\label{BM}
Let $a$ and $b$ be rational integers with $2\le b\le 10$ and $1\le a\le b-1$. The integer $N$ with all digits equal to $a$ in base $b$ is not a pure power, except for $N=1, 4, 8$ or $9$, for $N=11111$ written in base $b=3$, for $N=1111$ written in base $b=7$ and for $N=4444$ written in base $b=7$.
\end{thm}

\begin{thm}(Ljunggren \cite{Ljunggren})
\label{L}
The only integer solutions $(x,n,y)$ with $|x|>1$, $n>2$ and $y>0$ to the exponential equation 
\begin{align*}
\frac{x^n-1}{x-1}=y^2
\end{align*}
are $(x,n,y)=(7,4,20)$ and $(x,n,y)=(3,5,11)$.
\end{thm}

\section{Proofs}

\begin{proof}[Proof of Theorem 1.1]
a) Since any odd near-perfect number is square, the result follows from Theorem \ref{BMS}.

b) It easy to show that there are no near-perfect numbers of the form $p^k$, $k\ge0$ where $p$ is prime. Suppose that there exists an even near-perfect number of type A belonging to the Fibonacci sequence. For $p=2$ or $p=3$ one gets the numbers $18$ and $196$ which do not belong to the Fibonacci sequence.

Assume now that $p\ge5$. The equation $F_n=2^{p-1}(2^p-1)^2$ implies that $16|F_n$. From this it follows that $12|n$. Hence $3=F_4|F_n=2^{p-1}(2^p-1)^2$ which is impossible because $p\ge5$ and $2^p-1$ is prime. A similar argument can be used to show that there are no type B or type C near-perfect Fibonacci numbers.

Suppose now that $F_n$ is a near-perfect Fibonacci number with $\omega(F_n)=3$. Since $F_n$ is even, by Theorems \ref{P1} and \ref{P2} $n=3p$, $p>3$ and $F_n=2q_1q^{\alpha}_2$ where $F_p=q_1$ and $q_2$ is a primitive divisor of $F_n$, $\alpha\ge1$. If $q_1\ge7$, then 
\begin{align*}
2=\frac{\sigma(F_n)}{F_n}-\frac{d}{F_n}<\frac{3}{2}\cdot\frac{q_1+1}{q_1}\cdot\frac{q_2}{q_2-1}<\frac{3}{2}\cdot\frac{8}{7}\cdot\frac{11}{10}<2
\end{align*}
which is impossible. Thus, $q_1=5$. Then $F_n=F_{15}=2\cdot5\cdot61$ which is not a near-perfect number.

c) Clearly $L_6=18$ is a near-perfect number of type A. Using the fact that no member of the Lucas sequence is divisible by $8$ it is easy to verify that there are no other near-perfect Lucas numbers with two distinct prime divisors.
\end{proof}

\begin{proof}[Proof of Theorem 1.2]
For odd near-perfect numbers the result follows unconditionally from the Baker's Theorem \ref{B}. Note that if $m$ is sufficiently large near-perfect number with $\omega(m)=2$ then $\text{rad}(m)\ll\sqrt{m}$. Asssume $P(n)$ is a near-perfect number with large value of $n$, $\text{deg}P=d\ge3$ and $\omega(P(n))=2$. Fix $\epsilon>0$. Applying \eqref{Eq2.1}, we obtain 
\begin{align*}
n^{d-1-\epsilon}\ll\text{rad}(P(n))\ll n^{\frac{1}{2}d}.
\end{align*}
From the above, we get 
\begin{align*}
\frac{1}{2}d\ge d-1-\epsilon
\end{align*}
or $d\le2+\epsilon<3$. This contradiction implies the result.
\end{proof}

\begin{proof}[Proof of Theorem 1.3]
a) First we consider the near-perfect numbers of type A. Thus, to find repdigit near-perfect numbers we need to solve the equation 
\begin{align*}
N=d\bigg(\frac{g^n-1}{g-1}\bigg)=2^{p-1}(2^p-1)^2, \ \text{where} \ d\in\{1,\ldots,g-1\} \ \text{and} \ 2^p-1 \ \text{is prime}
\end{align*}
For the sake of contradiction assume that $n\ge3$. It's clear that $2^p-1\mid\frac{g^n-1}{g-1}$ for otherwise $(2^p-1)^2\mid d$ and then 
\begin{align*}
g>d\ge(2^p-1)^2>\sqrt{N}\ge\bigg(\frac{g^n-1}{g-1}\bigg)^{1/2}=\sqrt{g^{n-1}+\ldots+1}>g^{\frac{n-1}{2}}\ge g
\end{align*}
which is impossible. Thus $\frac{g^n-1}{g-1}=2^b(2^p-1)^2$ or $\frac{g^n-1}{g-1}=2^b(2^p-1)$ for some nonnegative integer $b$. Consider the first case. If $g$ is even, then $\frac{g^n-1}{g-1}$ is odd therefore $b=0$. Hence $\frac{g^n-1}{g-1}=(2^p-1)^2$ which has no solutions for $n\ge3$ by Theorem \ref{L}. Thus $g$ must be odd and $n$ must be even. Write $n=2m$. We then get
\begin{align*}
2^b(2^p-1)^2=\frac{g^{2m}-1}{g-1}=(g^m+1)\bigg(\frac{g^m-1}{g-1}\bigg)
\end{align*}
Note that $g^m+1>\frac{g^m-1}{g-1}$ and $2^p-1>2^b$. Moreover, $\text{gcd}\bigg(g^m+1,\frac{g^m-1}{g-1}\bigg)\le2$. Therefore, $g^m+1=2(2^p-1)^2$, $\frac{g^m-1}{g-1}=2^{b-1}$. The latter equation has no solutions in view of our assumption $2\le g\le10$ and Theorem \ref{BM}.

Now suppose that $\frac{g^n-1}{g-1}=2^b(2^p-1)$. If $g$ is even, then $\frac{g^n-1}{g-1}$ is odd therefore $b=0$. Hence
\begin{align*}
d=2^{p-1}(2^p-1)>2^p-1=\frac{g^n-1}{g-1}=g^{n-1}+\ldots+1>g^{n-1}>g
\end{align*}
which contradicts the assumption $1\le d\le g-1$. Thus $g$ must be odd and $n$ must be even. Put $n=2m$. We then obtain 
\begin{align*}
2^b(2^p-1)=\frac{g^{2m}-1}{g-1}=(g^m+1)\bigg(\frac{g^m-1}{g-1}\bigg)
\end{align*}
Since $g^m+1>\frac{g^m-1}{g-1}$ and $2^p-1>2^b$, $2^p-1\mid g^m+1$ and we get $g^m+1=2(2^p-1)$ and $\frac{g^m-1}{g-1}=2^{b-1}$. Since $\frac{g^m-1}{g-1}$ is even and $g$ is odd, we get that $m$ is even. Hence $m=2m_1$ and so $2(2^p-1)=g^m+1=g^{2m_1}+1\equiv2(\text{mod}\ 8)$. Then $2^p-1\equiv1(\text{mod} \ 4)$ but this is impossible for any prime $p\ge2$. Observe that for this case we did not use the assumption $2\le g\le10$. 

In a similar manner one can show finiteness of repdigits in base $g$ among near-perfect numbers of types B and C.

b) The result is an immediate consequence of Theorem \ref{BM}.
\end{proof}

Theorem \ref{Thm3} asserts that repdigit near-perfect numbers of types A, B and C have at most two digits in base $g$, $2\le g\le10$. For $g\in\{2,3,4,6\}$ there are no repdigit near-perfect numbers with two distinct prime factors. For $g=5$, 12, 18 and 24 are the only repdigit near-perfect numbers with two distinct prime factors. For $g=7$, 24 and 40 are the only repdigit near-perfect numbers with two distinct prime factors. For $g=8$, 18 is the only repdigit near-perfect number with two distinct prime factors. For $g=9$, 20 and 40 are the only repdigit near-perfect numbers with two distinct prime factors. Finally, in base $g=10$, 88 is the only repdigit near-perfect number with two distinct prime factors.



\normalsize

\begin{thebibliography}{10}

\bibitem{Baker}
A.~Baker,
\newblock{Bounds for the solutions of the hyperelliptic equation},
\newblock{\em Proc. Camb. Phil. Soc.}, {\bf 65} (1969), 439-444.

\bibitem{BGLLHT}  
K.~A. Broughan, M.~J. Gonzalez, R.~H. Lewis, F.~Luca, V.J.~M. Huguet and A.~Togbe,
\newblock{There are no multiply-perfect Fibonacci numbers},
\newblock{\em Integers}, {\bf 11} (2011), 363-397.

\bibitem{BM}
Y.~Bugeaud and M.~Mignotte,
\newblock{On integers with identical digits},
\newblock{\em Mathematika}, {\bf 46} (1999), no. 2, 411-417.

\bibitem{BMS}
Y.~Bugeaud, M.~Mignotte and S.~Siksek,
\newblock{Classical and modular approaches to exponential Diophantine equations I. Fibonacci and Lucas perfect powers},
\newblock{\em Ann. of Math.}, {\bf 163} (2006), no. 3, 969-1018.

\bibitem{Klurman}
D.~Klurman,
\newblock{Radical of perfect numbers and perfect numbers among polynomial values},
\newblock{\em Int. J. Number Theory}, {\bf 12} (2016), no. 3, 585-591.

\bibitem{Ljunggren}
W.~Ljunggren,
\newblock{Some theorems on indeterminate equations of the form $\frac{x^n-1}{x-1}=y^q$}
\newblock{\em Norsk Mat. Tidsskr.}, {\bf 25} (1943), 17-20.

\bibitem{Luca}
F.~Luca,
\newblock{Perfect Fibonacci and Lucas numbers}
\newblock{\em Rend. Circ. Mat. Palermo (2)}, {\bf 49} (2000), no. 2, 311-318.

\bibitem{LP}
F.~Luca and P.~Pollack,
\newblock{Multiperfect numbers with identical digits},
\newblock{\em J. Number Theory}, {\bf 131} (2011), no. 2, 260-284.


\bibitem{Pongsriiam}
P.~Pongsriiam,
\newblock{Fibonacci and Lucas numbers which have exactly three prime factors and some unique properties of $F_{18}$ and $L_{18}$},
\newblock{\em Fibo. Quart.}, {\bf 57} (2019), no. 5, 130-144.

\bibitem{Pollack}
P.~Pollack,
\newblock{Perfect numbers with identical digits},
\newblock{\em Integers}, {\bf 11} (2011), 519-529.

\bibitem{PS}
P.~Pollack and V.~Shevelev,
\newblock{On perfect and near-perfect numbers},
\newblock{\em J. Number Theory}, {\bf 132} (2012), 3037-3046.

\bibitem{RC}
X.~Z. Ren and Y.~G. Chen,
\newblock{On near-perfect numbers with two disctinct prime factors},
\newblock{\em Bull. Aust. Math. Soc.}, {\bf 88} (2013), 520-524.

\bibitem{TMF}
M.~Tang, X.~Y. Ma and M.~Feng
\newblock{On near-perfect numbers},
\newblock{\em Colloq. Math.}, {\bf 144} (2016), 157-188.

\bibitem{TRL}
M.~Tang, X.~Z. Ren and M.~Li
\newblock{On near-perfect and deficient perfect numbers},
\newblock{\em Colloq. Math.}, {\bf 133} (2013), 221-226.


\end{thebibliography}
\end{document}